\documentclass{amsart}
\usepackage{amssymb,amsthm,amsmath,epsfig,latexsym,xypic}
\usepackage{calc}
\usepackage{latexsym}
\usepackage{amscd,amssymb,subfigure,hyperref}
\usepackage[arrow,matrix,graph,frame,poly,arc,tips]{xy}

\begin{document}

\newcommand{\mmbox}[1]{\mbox{${#1}$}}
\newcommand{\affine}[1]{\mmbox{{\mathbb A}^{#1}}}
\newcommand{\Ann}[1]{\mmbox{{\rm Ann}({#1})}}
\newcommand{\caps}[3]{\mmbox{{#1}_{#2} \cap \ldots \cap {#1}_{#3}}}
\newcommand{\N}{{\mathbb N}}
\newcommand{\Z}{{\mathbb Z}}
\newcommand{\R}{{\mathbb R}}
\newcommand{\KK}{{\mathbb K}}
\newcommand{\A}{{\mathcal A}}
\newcommand{\B}{{\mathcal B}}
\newcommand{\OO}{{\mathcal O}}
\newcommand{\C}{{\mathbb C}}
\newcommand{\PP}{{\mathbb P}}

\newcommand{\Tor}{\mathop{\rm Tor}\nolimits}
\newcommand{\Ext}{\mathop{\rm Ext}\nolimits}
\newcommand{\Hom}{\mathop{\rm Hom}\nolimits}
\newcommand{\im}{\mathop{\rm im}\nolimits}
\newcommand{\rk}{\mathop{\rm rk}\nolimits}
\newcommand{\codim}{\mathop{\rm codim}\nolimits}
\newcommand{\supp}{\mathop{\rm supp}\nolimits}
\newcommand{\coker}{\mathop{\rm coker}\nolimits}
\sloppy
\newtheorem{defn0}{Definition}[section]
\newtheorem{prop0}[defn0]{Proposition}
\newtheorem{conj0}[defn0]{Conjecture}
\newtheorem{thm0}[defn0]{Theorem}
\newtheorem{lem0}[defn0]{Lemma}
\newtheorem{corollary0}[defn0]{Corollary}
\newtheorem{example0}[defn0]{Example}

\newenvironment{defn}{\begin{defn0}}{\end{defn0}}
\newenvironment{prop}{\begin{prop0}}{\end{prop0}}
\newenvironment{conj}{\begin{conj0}}{\end{conj0}}
\newenvironment{thm}{\begin{thm0}}{\end{thm0}}
\newenvironment{lem}{\begin{lem0}}{\end{lem0}}
\newenvironment{cor}{\begin{corollary0}}{\end{corollary0}}
\newenvironment{exm}{\begin{example0}\rm}{\end{example0}}

\newcommand{\msp}{\renewcommand{\arraystretch}{.5}}
\newcommand{\rsp}{\renewcommand{\arraystretch}{1}}

\newenvironment{lmatrix}{\renewcommand{\arraystretch}{.5}\small
 \begin{pmatrix}} {\end{pmatrix}\renewcommand{\arraystretch}{1}}
\newenvironment{llmatrix}{\renewcommand{\arraystretch}{.5}\scriptsize
 \begin{pmatrix}} {\end{pmatrix}\renewcommand{\arraystretch}{1}}
\newenvironment{larray}{\renewcommand{\arraystretch}{.5}\begin{array}}
 {\end{array}\renewcommand{\arraystretch}{1}}

\title
{The weak Lefschetz property and powers of linear forms in $\KK[x,y,z]$}

\author{Hal Schenck}
\thanks{Schenck supported by NSF 07--07667, NSA 904-03-1-0006}
\address{Schenck: Mathematics Department \\ University of
 Illinois \\
   Urbana \\ IL 61801\\USA}
\email{schenck@math.uiuc.edu}

\author{Alexandra Seceleanu}
\address{Seceleanu: Mathematics Department \\ University of
 Illinois \\
   Urbana \\ IL 61801\\USA}
\email{asecele2@math.uiuc.edu}


\begin{abstract}
\noindent We show that an Artinian quotient of an ideal $I \subseteq \KK[x,y,z]$
 generated by powers of linear forms has the Weak Lefschetz property. If the syzygy bundle 
of $I$ is semistable this follows from results of Brenner-Kaid; our proof works without
 this hypothesis, which typically does not hold.
\end{abstract}
\maketitle
 
\section{Introduction}\label{sec:one}
Let $I \subseteq S=\KK[x_1,\ldots, x_r]$ be an ideal such that $A=S/I$ is Artinian. Then $A$ has 
the {\em Weak Lefschetz Property} (WLP) if there is an $\ell \in S_1$ such that for all $m$, the map $\mu_{\ell}$ 
\[
A_m \stackrel {\cdot \ell}{\longrightarrow}A_{m+1}
\]
is either injective or surjective. We assume $char(\KK)=0$; as shown in \cite{MMN}, WLP behaves
in very subtle ways in positive characteristic. In \cite{A}, 
Anick shows that if $r=3$ and $I$ is generated by 
generic forms, then $A$ has WLP. In \cite{HMNW}, Harima-Migliore-Nagel-Watanabe introduced
the syzygy bundle of $I$ to study the WLP, and this bundle also plays a key role  in recent work of Brenner-Kaid \cite{BK}.
\begin{defn}
If $I = \langle f_1,\ldots, f_n \rangle$ is $\langle x_1,\ldots,x_r \rangle-$primary, and $deg(f_i) = d_i$, 
then the syzygy bundle ${\mathcal S}(I) = \widetilde{Syz(I)}$ is a rank $n-1$ bundle defined via:
\[
0 \longrightarrow Syz(I) \longrightarrow \bigoplus\limits_{i=1}^{n}S(-d_i) \stackrel{[f_1,\ldots,f_n]}{\longrightarrow} S
\]
\end{defn}
The cokernel of the rightmost map is $S/I$, which vanishes as a sheaf.
\begin{defn}
A vector bundle $\mathcal{E}$ on projective space is said to be semistable if
for every coherent subsheaf $\mathcal{F} \subseteq \mathcal{E}$
\[
\frac{c_1(\mathcal{F})}{rk(\mathcal{F})} \le \frac{c_1(\mathcal{E})}{rk(\mathcal{E})}, \mbox{ where }c_1\mbox{ denotes the first Chern class}.
\]
\end{defn}
By a result of Grothendieck, every vector bundle on $\PP^1$ splits as a sum of line bundles \cite{OK}, so for a given line $L$, if $\mathcal{E}$ 
has rank $k$, then
\[
\mathcal{E}|_L \simeq \bigoplus\mathcal{O}_L(a_i), \mbox{ with } a_1 \ge a_2 \ge \cdots \ge a_k.
\]

If $\mathcal{E}$ is semistable, then \cite{OK} for a generic line $L$ the 
tuple $(a_1, a_2, \cdots, a_k)$ does not vary;  $(a_1, a_2, \cdots, a_k)$ 
is the {\em generic splitting type} of $\mathcal{E}$; if $\mathcal{E}$ is 
semistable, then $|a_i-a_{i+1}| \le 1$.
\pagebreak

For the remainder of the paper we focus on the case $r=3$, so henceforth
$S$ denotes $\KK[x,y,z]$. In \cite{BK}, Brenner and Kaid show 
that if $A = S/I$ is Artinian with ${\mathcal S}(I)$ semistable of generic splitting type
$(a_1,\ldots, a_{n-1})$, then $A$ has WLP iff $|a_1 -a_{n-1}|\le 1$. As a 
corollary of this, they recover a result of Harima-Migliore-Nagel-Watanabe \cite{HMNW}
that every Artinian complete intersection in $S$ has WLP. They also 
completely characterize WLP for almost complete intersections, showing that
in this case if ${\mathcal S}(I)$ is not semistable, then WLP holds.

It is clear from the definition that semistability can be a difficult property to
show. In this note, we examine a special class of ideals in $S$ which falls outside
the classes considered by Anick, Brenner-Kaid, and Harima-Migliore-Nagel-Watanabe.
Our main result is \vskip .05in
\noindent{\bf Theorem} An Artinian quotient of $\KK[x,y,z]$ by powers of linear forms has WLP.

\section{Proof of the theorem}\label{sec:two}
We begin by recalling the setup of \cite{BK}. Let $\ell$ be a generic form in $S_1$ with
$L = V(\ell)$, and $I$ an ideal such that $A=S/I$ is Artinian. Taking cohomology of the defining sequence for ${\mathcal S}(I)$ 
\[
0 \longrightarrow {\mathcal S}(I)(m) \longrightarrow \bigoplus\limits_{i=1}^{n}\OO_{\PP^2}(m-d_i) \longrightarrow \OO_{\PP^2}(m) \longrightarrow 0,
\]
we see that 
\[
A = \bigoplus_{m \in \Z} H^1({\mathcal S}(I)(m)).
\]
On the other hand, since ${\mathcal S}(I)$ is a bundle, tensoring the
sequence 
\[
0 \longrightarrow \OO_{\PP^2}(m) \longrightarrow \OO_{\PP^2}(m+1) \longrightarrow \OO_L(m+1) \longrightarrow 0
\]
with ${\mathcal S}(I)$ gives the exact sequence
\[
0 \longrightarrow {\mathcal S}(I)(m) \longrightarrow {\mathcal S}(I)(m+1) \longrightarrow {\mathcal S}(I)|_L(m+1) \longrightarrow 0
\]
The long exact sequence in cohomology yields a sequence

\[
\xymatrixrowsep{25pt}
\xymatrixcolsep{25pt}
\xymatrix{
0 \ar[r] & H^0({\mathcal S}(I)(m)) \ar[r] &H^0({\mathcal S}(I)(m+1)) \ar[r]^{\phi_m} &H^0({\mathcal S}(I)|_L(m+1)) \ar[dll]\\
         & H^1({\mathcal S}(I)(m)) \ar[r]^{\cdot \ell} &H^1({\mathcal S}(I)(m+1)) \ar[r] &H^1({\mathcal S}(I)|_L(m+1)) \ar[dll]^{\psi_m}\\
         & H^2({\mathcal S}(I)(m)) \ar[r] &H^2({\mathcal S}(I)(m+1)) \ar[r] &H^2({\mathcal S}(I)|_L(m+1)) = 0.
}
\]
Therefore injectivity of $\mu_{\ell}$ follows from surjectivity of $\phi_m$, and
surjectivity of $\mu_{\ell}$ from injectivity of $\psi_m$. Our next step is to analyze ${\mathcal S}(I)|_L$. To do this,
we tensor the defining sequence 
\[
0 \longrightarrow Syz(I) \longrightarrow \bigoplus\limits_{i=1}^n S(-d_i) \longrightarrow I \longrightarrow 0
\] 
with $S/\ell$, yielding the sequence
\[
0 \longrightarrow Tor^S_1(I,S/\ell)\longrightarrow  Syz(I)\otimes S/\ell \longrightarrow \bigoplus\limits_{i=1}^n S/\ell(-d_i) \longrightarrow I\otimes S/\ell \longrightarrow 0.
\] 
Now $Tor^S_1(I,S/\ell) = 0$ since it is the kernel of 
\[
I \stackrel{\cdot \ell}{\longrightarrow} I(1).
\]
After a change of coordinates, $\ell = x$ is generic. Reducing the defining equations 
of $I$ mod $x$, we see that $Syz(I)\otimes S/\ell$ is the module of syzygies on 
$I\otimes S/\ell$, an ideal generated 
by powers of linear forms in two variables. We make use of the following pair of 
lemmas from \cite{GS} on ideals
\[
J = \langle l_1^{\alpha _1},\ldots, l_t^{\alpha_t}\rangle \subseteq \KK[y,z]=R,
\]
generated by powers of pairwise linearly independent forms.
\vskip .05in
\begin{lem}~\label{mingens}
Let $0 < \alpha _1 \leq \alpha _2 \cdots \leq \alpha _t$. Then for $m \ge 2$:
\[
l_{m+1}^{\alpha _{m+1}} \notin \langle l_1^{\alpha _1}, \ldots , l_m^{\alpha _m}\rangle
\Leftrightarrow \alpha_{m+1} \le \frac{\sum _{i=1}^m \alpha_i - m}{m-1}.
\]
\end{lem}
\vskip .05in
\begin{lem}~\label{syzs}
If $J$ is minimally generated by $\{l_1^{\alpha_1}, \ldots , l_t^{\alpha_t}\}$ and $t\ge 2$, 
then the socle degree of $\KK[y,z]/J$ is 
$\omega = \Big\lfloor \frac{\sum_{i=1}^t\alpha_i -t}{t-1}\Big\rfloor$, and $J$ has 
minimal free resolution 
\[
0 \longrightarrow  R(-\omega-2)^a \oplus R(-\omega-1)^{t-1-a} \longrightarrow \oplus_{i=1}^t R(-\alpha_i) \longrightarrow J \longrightarrow 0,
\]
where 
$ a = \sum_{i=1}^t\alpha_i -(t-1)(\omega+1).$
\end{lem}
\vskip .05in
\begin{prop}~\label{noredun}
If $I = \langle l_1^{d_1},\ldots,l_n^{d_n}\rangle \subseteq S$ satisfies 
\[
d_{t+1} \le \frac{\sum _{i=1}^t d_i - t}{t-1},
\]
for all $t > 1$, then $S/I$ has the WLP.
\end{prop}
\begin{proof}
By Lemma~\ref{mingens}, the restriction $I \otimes S/\ell$ has the same number of minimal
generators and degrees as $I$, and so it follows from Lemma~\ref{syzs} that
\begin{equation}\label{Split}
{\mathcal S}(I)|_L \simeq \OO_L(-\omega-2)^a \oplus \OO_L(-\omega-1)^{n-1-a},
\end{equation}
with 
\[\omega = \Big\lfloor \frac{\sum_{i=1}^n d_i -n}{n-1}\Big\rfloor \mbox{  and }a = \sum_{i=1}^n d_i -(n-1)(\omega+1).
\]
Suppose $m < \omega$. Then 
\[
H^0({\mathcal S}(I)|_L(m+1) \simeq H^0(\OO_L(m-1-\omega))^a \oplus H^0(\OO_L(m-\omega))^{n-1-a} = 0,
\]
so $\mu_{\ell}$ is injective. If instead $m \ge \omega$, by Serre duality
\[
H^1({\mathcal S}(I)|_L(m+1) \simeq H^0(\OO_L(-m-1+\omega))^a \oplus H^0(\OO_L(-m-2+\omega))^{n-1-a} = 0,
\]
and thus $\mu_{\ell}$ is surjective.
\end{proof}
\pagebreak
\begin{thm}~\label{mainT}
If $I = \langle l_1^{d_1},\ldots,l_n^{d_n}\rangle \subseteq S$, then $S/I$ has the WLP.
\end{thm}
\begin{proof}
If 
\[
d_{t+1} \le \frac{\sum _{i=1}^t d_i - t}{t-1},
\]
for all $t > 1$, then this follows from Proposition~\ref{noredun}, so let $d_1 \le d_2\le \cdots \le d_n$ and suppose that $t+1$ is the first index where 
\[
d_{t+1} > \frac{\sum _{i=1}^t d_i - t}{t-1} \ge \Big\lfloor \frac{\sum_{i=1}^t d_i -t}{t-1}\Big\rfloor = \omega.  
\]
Thus, $d_i \ge \omega +1$ when $i \ge t+1$. If $d_i \in \{\omega +1, \omega +2 \}$ for all
$i \ge t+1$, then the shifts appearing in ${\mathcal S}(I)|_L$ are as in Equation~\ref{Split},
so the argument of Proposition~\ref{noredun} works. Suppose $k \ge t+1$ is the first
index such that $d_k \ge \omega +3$. Then 
\[
{\mathcal S}(I)|_L \simeq \OO_L(-\omega-1)^{a} \bigoplus \OO_L(-\omega-2)^{b} \bigoplus\limits_{i=k}^n\OO_L(-d_i),
\]
with $a+b =k-2$. If $m < \omega$, the argument of Proposition~\ref{noredun}
shows that $\mu_{\ell}$ is injective, so suppose $m \ge \omega$. We show $\psi_m$ is injective by
a dimension computation. From the defining sequence for ${\mathcal S}(I)$ we obtain
\[
0 \longrightarrow H^2({\mathcal S}(I)(m)) \longrightarrow H^2(\bigoplus\limits_{i=1}^{n}\OO_{\PP^2}(m-d_i)) \longrightarrow  H^2(\OO_{\PP^2}(m)) \longrightarrow 0.
\]
By Serre duality, $h^2(\OO_{\PP^2}(m)) = h^0(\OO_{\PP^2}(-m-3))= 0$ since $m \ge \omega >0$, and   
\[
h^2({\mathcal S}(I)(m)) = \sum\limits_{i=1}^n {d_i-m-1 \choose 2} \mbox{ and } h^2({\mathcal S}(I)(m+1)) = \sum\limits_{i=1}^n {d_i-m-2 \choose 2}.
\]
Thus, 
\[
\dim \im (\psi_m) = \sum\limits_{i=1}^n \max(d_i-m-2, 0).
\]
The contributions come from those $d_i \ge m+3 \ge \omega+3$. Our assumption is that
\[
{\mathcal S}(I)|_L \simeq \OO_L(-\omega-1)^{a} \bigoplus \OO_L(-\omega-2)^{b} \bigoplus\limits_{i=k}^n\OO_L(-d_i),
\]
with $a+b =k-2$. Thus for $m \ge \omega$, 
\[
\begin{array}{cc}
\!\!\!h^1({\mathcal S}(I)|_L(m\!+\!1) &\!\!\! =\sum\limits_{i=k}^n h^1(\OO_L(-d_i\!+\!m\!+\!1))\!+\! h^1(\OO_L(m\!-\!\omega)^{a})\!+\! h^1(\OO_L(m\!-\!\omega\!-\!1)^{b})\\
                   &\!\!\!\!\!\!\!\!\!\!\!\!\!\!\!\!\!\!\!\!\!\!\!\!\!\!\!\!\!\!\!\!\!\!\!\!\!\!\!\!\!\!\!\!\!\!\!\!\!\!\!\!\!\!\!\!\!\!\!\!\!\!\!\!\!\!\!\!\!\!\!\!\!\!\!\!\!\!\!\!\!\!\!\!\!\!\!\!\!\!\!\!\!\!\!\!\!\!\!\!\!\!\!=\sum\limits_{i=k}^n h^1(\OO_L(-d_i\!+\!m\!+\!1))\\
                   &\!\! \!\!\!\!\!\!\!\!\!\!\!\!\!\!\!\!\!\!\!\!\!\!\!\!\!\!\!\!\!\!\!\!\!\!\!\!\!\!\!\!\!\!\!\!\!\!\!\!\!\!\!\!\!\!\!\!\!\!\!\!\!\!\!\!\!\!\!\!\!\!\!\!\!\!\!\!\!\!\!\!\!\!\!\!\!\!\!\!\!\!\!\!\!\!\!\!\!\!\!\!\!\!= \sum\limits_{i=k}^n h^0(\OO_L(d_i-m-3))  \\
                   &\!\!\!\!\!\!\!\!\!\!\!\!\!\!\!\!\!\!\!\!\!\!\!\!\!\!\!\!\!\!\!\!\!\!\!\!\!\!\!\!\!\!\!\!\!\!\!\!\!\!\!\!\!\!\!\!\!\!\!\!\!\!\!\!\!\!\!\!\!\!\!\!\!\!\!\!\!\!\!\!\!\!\!\!\!\!\!\!\!\!\!\!\!\!\!\!\!\!\!\!\!\!\!\!\!= \sum\limits_{i=k}^n \max(d_i-m-2, 0).
\end{array}
\]
Since this is equal to $\dim \im (\psi_m)$, $\psi_m$ is an inclusion, so that $\mu_{\ell}$ is
surjective.
\end{proof}
\pagebreak
It follows from Theorem~\ref{mainT} that ideals generated by 
powers of linear forms in $\KK[x,y,z]$ which have unstable syzygy bundles 
always have WLP. As noted earlier Brenner and Kaid show that almost complete
intersections with unstable syzygy bundles have WLP. Thus, it seems
reasonable to ask if every ideal in $\KK[x,y,z]$ with unstable 
syzygy bundle has WLP. 
\begin{exm}\label{nonSSnonWLP}
For the ideal $I = \langle x^5,y^5,z^5,x^2yz,xy^2z \rangle \subseteq \KK[x,y,z]$, 
${\mathcal S}(I)$ is not semistable, by Proposition~2.2 of \cite{B}. The Hilbert 
function of $A$ is $(1,3,6,10,13,13,10,6,3)$ and a computation shows the map 
from $A_4 \rightarrow A_5$ is not full rank, so $A$ does not have WLP.  
\end{exm}
\noindent As noted, Theorem~\ref{mainT} need not hold for more than three variables:
\begin{exm}\label{nonWLPlinforms}
The ring $A = \KK[x,y,z,w]/\langle x^3,y^3,z^3,w^3,(x+y+z+w)^3\rangle$ appears
in Example~8.1 of \cite{MMN}, and does not have WLP. The Hilbert function of $A$ is $(1,4,10,15,15,6)$, and a computation shows the map from $A_3 \rightarrow A_4$ is not full rank. So WLP need not hold for powers of 
linear forms in more than three variables. 
\end{exm}
\vskip .05in
\noindent{\bf Concluding Remarks} $\mbox{ }$The proof of Theorem~\ref{mainT} works for any ideal which has the
same splitting type as an ideal generated by powers of linear forms, so
it would be interesting to find families of such ideals. In light 
of Example~\ref{nonWLPlinforms}, we ask: are there reasonable additional 
hypotheses so that a version of 
Theorem~\ref{mainT} does hold in more than three variables? A second
question is if ideals generated by powers of linear forms 
possess the Strong Lefschetz Property. As pointed out by the referee,
the answer is no: SLP fails for the ideal generated by cubes of four
general linear forms, and multiplication by a cube of a linear form.
However, multiplication by a general form of degree three {\em does}
have maximal rank, so we ask: does multiplication by a general form
of any degree induce a multiplication having maximal rank?

\vskip .05in
\noindent{\bf Acknowledgements} Computations were performed using Macaulay2,
by Grayson and Stillman, available at: {\tt http://www.math.uiuc.edu/Macaulay2/}.
Scripts to analyze WLP are available at: {\tt http://www.math.uiuc.edu/$\sim$asecele2}.
We thank an anonymous referee for thoughtful comments.
\bibliographystyle{amsalpha}

\begin{thebibliography}{10}

\bibitem{A} D. Anick, 
            Thin algebras of embedding dimension three,
            {\em J. Algebra}, {\bf 100} (1986), 235--259.

\bibitem{B} H.~Brenner,
            Looking out for stable syzygy bundles. 
            {\em Adv. Math.} {\bf 219}  (2008), 401--427.

\bibitem{BK} H.~Brenner, A.~Kaid, 
            Syzygy bundles on $\mathbb{P}^2$ and the weak Lefschetz property.  
            {\em Illinois J. Math.} {\bf 51}  (2007), 1299--1308.

\bibitem{GS} A.~Geramita, H.~Schenck, 
            Fat points, inverse systems, and piecewise polynomial functions,
            {\em J. Algebra}, {\bf 204} (1998), 116--128.

\bibitem{HMNW} T.~Harima, J.~Migliore, U.~Nagel, J.~Watanabe,
            The weak and strong Lefschetz properties for Artinian $\KK$-algebras.  
            {\em J. Algebra}, {\bf 262} (2003), 99--126.

\bibitem{MMN} J.~Migliore, R.~Mir\'o-Roig,  U.~Nagel,
            Monomial ideals, almost complete intersections and the weak Lefschetz property
            {\em Trans. Amer. Math. Soc.},to appear.

\bibitem{OK} C.~Okonek, M.~Schneider, H.~Spindler,
        {\em Vector Bundles on Complex Projective Spaces}, Progress in Mathematics,vol. ~3,
        Birkhauser, Boston, 1980.

\end{thebibliography}

\end{document}